\documentclass{amsart}
\usepackage{graphicx}
\newtheorem{thm}{Theorem}[section]

\newtheorem{prop}[thm]{Proposition}
\theoremstyle{definition}

\theoremstyle{remark}

\numberwithin{equation}{section}

\begin{document}

\title[]{ON SINGULAR EXTENSIONS OF CONTINUOUS FUNCTIONALS FROM C([0,1])
TO VARIABLE LEBESGUE SPACES}%

\author{Daviti Adamadze}
\address{Daviti Adamadze \\
Faculty of Exact and Natural Sciences\\
 Javakhishvili Tbilisi State
University\\
 13, University St., Tbilisi, 0143, Georgia}
\email{daviti.adamadze2013@ens.tsu.edu.ge}

\author{Tengiz Kopaliani}
\address{Tengiz Kopaliani \\
Faculty of Exact and Natural Sciences\\
 Javakhishvili Tbilisi State
University\\
 13, University St., Tbilisi, 0143, Georgia}
\email{tengiz.kopaliani@tsu.ge}


\thanks{The research of the second  author was  supported by Shota Rustaveli National Science Foundation of Georgia
(SRNSFG) FR$17_{-}589$.}

\begin{abstract}
  Valadier and  Hensgen proved independently that the restriction of functional $\phi(x)=\int_{0}^{1}x(t)dt,\,\,x\in L^{\infty}([0,1])$ on the space of  continuous functions $C([0,1])$ admits a singular extension back to the whole space $L^{\infty}([0,1]).$  Some general results in this direction for the Banach lattices were obtained by Abramovich and Wickstead. In present note we investigate analogous problem for variable exponent Lebesgue spaces, namely we prove that if  the space of continuous functions $C([0,1])$ is closed subspace in  $L^{p(\cdot)}([0,1]),$ then every bounded linear  functional on $C([0,1])$ is the restriction of a singular linear functional on $L^{p(\cdot)}([0,1])$.

\end{abstract}
\maketitle
\section{Introduction}

In the theory of duality of function spaces an investigation of the space of all
singular linear functionals is of importance. It is well known that the topological dual of the Banach  function space $X$ can be represented in the form
$X^{\ast}=X_{n}^{\ast}\oplus X_{s}^{\ast}$, where $X_{n}^{\ast}$ is the order continuous dual of $X$  and $X_{s}^{\ast}$ (the disjoint complement of $X_{n}^{\ast}$ in $X^{\ast})$ is the space of all singular linear functionals on $X.$

The term "singular functional" is overused in the literature. According to \cite{AB} the space of singular functionals $X_{s}^{\ast},$ is defined as we mentioned above as the band
$(X_{n}^{\ast})^{d},$ complementary to the band of order continuous  functionals.  Often the space of singular functionals for Banach function space $X$ is defined as an annihilator
$$
(X_{a})^{\perp}=\{x^{\ast}\in X^{\ast};\,\,x^{\ast}(x)=0,\,\,\,\mbox{for all}\,\,x\in X_{a}\},
$$
where $X_{a}$  is the space of order continuous elements in $X.$  Although, the differences between these definitions are small and they often define the same  objects (see, for instance \cite{AW}), for variable exponent Lebesgue spaces mentioned definitions are equivalent. (Note that for $L^{\infty}([0,1])$, we have $(L^{\infty}([0,1]))_{a}=\{0\}).$

Let  $L^{\infty}([0,1])$ and $C([0,1])$ denote respectively the Banach space of essentially bounded real-valued functions  on the $[0,1]$  with the Lebesgue measure, and its subspace of continuous functions.  Valadier \cite{Val} and  Hensgen \cite{He} proved independently  that the restriction of the functional
 $$
\phi(x)=\int_{0}^{1}x(t)dt,\,\,\,\,x\in L^{\infty}[0,1]
$$
to a fairly large subspace $C([0,1])$ of  $L^{\infty}([0,1])$  admits a singular ("bad") extension back to the $L^{\infty}([0,1])$. Abramovich and Wickstead showed, that every bounded linear functional on $C([0,1])$ is the restriction of a singular functional on $L^{\infty}([0,1]).$ They also generalized this result to the Banach lattice setting (\cite{AW}, see Theorem 1 and remarks thereafter). Let a finitely additive measure $\nu$ represents $f\in (L^{\infty}([0,1])^{\ast}$ and $\widehat{\nu}$ is the Borel measure representing $f$ restricted on  $C([0,1]).$   Many properties of $\widehat{\nu}$  in terms $\nu$  recently was inestigated by Toland \cite{Tol} and Wrobel \cite{WR}.

 Edmunds, Gogatishvili and Kopaliani \cite{EdGoKo} showed  that there is a variable
exponent space $L^{p(\cdot)}([0,1])$  with $1<p(t)<\infty$ a.e.,  which has in common with
$L^{\infty}([0,1])$ the property that the space $C([0,1])$ is
closed linear subspace in it. Moreover, both the Kolmogorov and the Marcinkiewicz
examples of functions with a.e. divergence Fourier series belong to $L^{p'(\cdot)}([0,1])$,
where $p'(\cdot)$ is conjugate function of $p(\cdot)$.

 In \cite{KoZv} there is given a
necessary and sufficient condition on the decreasing rearrangement $p^{\ast}$ of exponent
$p(\cdot)$, for existence of equimeasurable exponent function of $p(\cdot)$ which corresponding
variable Lebesgue space has the property that the space of continuous functions is
closed in it. Indeed, let for functions $p(\cdot):[0,1]\rightarrow[1,\infty)$ we have
$$
\lim\mbox{sup}_{t\rightarrow0_{+}}\frac{p^{\ast}(t)}{\ln(e/t)}>0,
$$
then there exists equimeasurable with $p(\cdot)$ exponent function $\widetilde{p}(\cdot)$  such that the space $C([0,1])$ is closed subspace in $L^{\widetilde{p}(\cdot)}([0,1])$.

Let the space $C([0,1])$ is closed subspace of $L^{p(\cdot)}([0,1])$. It is interesting to  investigate of validity of  analogous  theorem of Abramovich and Wickstead which was mentioned above, in this case. We gave answer about this question.We prove following
\begin{thm}
Let the space $C([0,1])$ is closed subspace of $L^{p(\cdot)}([0,1])$. Then every bounded linear functional on $C([0,1])$ is the restriction of a linear singular functional on $L^{p(\cdot)}([0,1])$.
\end{thm}

\section{ some properties of singular functionals in variable lebesgue spaces }

Let $p(\cdot):[0,1]\rightarrow [1,\infty)$ be a  measurable function. Define the modular
$$
\rho_{p(\cdot)}(x)=\int_{[0,1]}|x(t)|^{p(t)}dt.
$$
  Given a measurable function $x,$ we say that $x\in L^{p(\cdot)}([0,1])$ if there exists $\lambda>0$ such that $\rho_{p(\cdot)}(x/\lambda)<\infty.$
This set becomes a Banach function space when equipped with the Luxemburg  norm
$$
\|x\|_{p(\cdot)}=\inf\{\lambda>0;\,\,\,\rho_{p(\cdot)}(x/\lambda)\leq1\}.
$$
The variable Lebesgue spaces were first  introduced by Orlicz. They have been widely studied for the past thirty years, both for their interest as function spaces and for their applications to PDEs and the calculus of variation (see \cite{CUF},\cite{DHHR}).

 Define the dual exponent $p'(\cdot)$ pointwise by $1/p(t)+1/p'(t)=1,\,\,t\in[0,1].$

In the case $p_{+}<\infty$, where $
p_{+}=\mbox{ess\,sup}_{t\in [0,1]}\,p(t)$,  the dual space of $L^{p(\cdot)}([0,1])$ can be completely characterized, it is isomorphic to $L^{p'(\cdot)}([0,1]).$ The problem of characterizing the dual of $L^{p(\cdot)}([0,1])$ when $p_{+}=\infty$ was considered in \cite{ACCUO}. The authors in this case give a decomposition of $(L^{p(\cdot)}([0,1]))^{\ast}$ as a direct sum of $L^{p'(\cdot)}([0,1])$ and the dual of a quotient space (the authors refer to as the germ space, and denote by $L^{p(\cdot)}_{\mbox{germ}}$).  Note that some main aspect of this subject was done in more general setting for Musielak-Orlicz spaces by Hudzik and Zbaszyniak in the paper \cite{Hud}.  At first we present some basic facts from mentioned paper  for variable Lebesgue setting.  We will always assume without loss of generality that $1< p(t)<\infty$ a.e. (we are interested in characterizing the spaces $L^{p(\cdot)}([0,1])$ close to $L^{\infty}([0,1])$).

We define the closed subspace $E^{p(\cdot)}([0,1])$  of $L^{p(\cdot)}([0,1])$ by following
$$
E^{p(\cdot)}([0,1])=\{x:\,\,\rho_{p(\cdot)}(\lambda x)<\infty\,\,\,\,\mbox{for\,any}\,\,\lambda>0\}.
$$

It is easy to see that $E^{p(\cdot)}([0,1])$ is the subspace of order continuous elements in $L^{p(\cdot)}([0,1])$, i.e. $x\in L^{p(\cdot)}([0,1])$ belongs to $E^{p(\cdot)}([0,1])$ if and only if for any sequence $x_{n}$ of measurable functions on $[0,1]$ such that $|x_{n}(t)|\leq |x(t)|$ for all $n\in \mathbb{N}$  and $|x_{n}|\rightarrow 0$ a.e. on $[0,1]$ there holds $\|x_{n}\|_{p(\cdot)}\rightarrow0.$ For the definition of order continuous elements in Banach lattices see \cite{AB}.  Note that if $p_{+}<\infty$ the spaces $L^{p(\cdot}([0,1])$ and $E^{p(\cdot}([0,1])$ coincide each other (see \cite{CUF}, \cite{DHHR}).

Let $p_{+}=\infty.$ Define the sets $\Omega_{n}=\{t\in [0,1]:\,\,p(t)\leq n\},\,\,n\in \mathbb{N}.$ We will
always assume without loss of generality that $|\Omega_{n}|>0$ for $n\in \mathbb{N},\,\,n\geq 2$ and $|\Omega_{1}|=0.$ For $x\in L^{p(\cdot)}([0,1])$ define the  functions $x^{(n)}\in E^{p(\cdot)}([0,1]),\,\,n\in \mathbb{N}$ as $x^{(n)}=x\chi_{\Omega_{n}}$ ($\chi_{\Omega_{n}}$ denotes the characteristic function of the set $\Omega_{n}$).

For any $x\in L^{p(\cdot)}([0,1])$ define
$$
d(x)=\inf\{\|x-y\|_{p(\cdot)}:\,\,y\in E^{p(\cdot)}\},
$$

$$
\theta(x)=\inf\{\lambda>0;\,\,\rho_{p(\cdot)}(\lambda x)<+\infty\}.
$$

For any $x^{\ast}\in(L^{p(\cdot)}([0,1]))^{\ast}$ we define  the  norm in a dual space

$$
\|x^{\ast}\|=\sup\{x^{\ast}(x):\,\,\|x\|_{p(\cdot)}\leq1\}.
$$

The dual space of $L^{p(\cdot)}([0,1])$ is represented in the following way (see \cite{Hud}):
$$
(L^{p(\cdot)}([0,1]))^{\ast}=L^{p'(\cdot)}([0,1])\oplus(L^{p(\cdot)}([0,1]))^{\ast}_{s},
$$
i.e. every $x^{\ast}\in (L^{p(\cdot)}([0,1]))^{\ast}$ is uniquely represented in the form $x^{\ast}=\xi_{v}+\varphi$, where $\xi_{v}$ is the regular functional defined by a function  $v\in L^{p'(\cdot)}([0,1])$ by the formula
\begin{equation}
\xi_{v}(x)=\int_{[0,1]}v(t)x(t)dt,\,\,\,x\in L^{p(\cdot)}([0,1]),
\end{equation}
and $\varphi$ is a singular functional, i.e. $\varphi(x)=0$ for any $x\in E^{p(\cdot)}([0,1])$ (for the case $p_{+}<\infty$, we have $ L^{p(\cdot)}[0,1])^{\ast}_{s}=\{0\}$).

\begin{prop} (\cite{Hud}, Lemma 1.2)
For any $x\in L^{p(\cdot)}([0,1])$  there holds the equalities
$$
\lim_{n\rightarrow\infty}\|x-x_{n}\|_{p(\cdot)}=\theta(x)=d(x).
$$
\end{prop}

\begin{prop} (\cite{Hud}, Lemma 1.3)  For any  singular functional $\varphi$ there holds
$$
\|\varphi\|=\sup\{\varphi(x):\,\,\rho_{p(\cdot)}(x)<\infty\}=\sup_{x\in L^{p(\cdot)}\backslash E^{p(\cdot)}}\varphi(x)/\theta(x).
$$
\end{prop}

\begin{prop} (\cite{Hud}, Lemma 1.4)

 For any functional $x^{\ast}=\xi_{v}+\varphi\in (L^{p(\cdot)}([0,1]))^{\ast}$, where $\xi_{v}$ is defined by (2.1) and $\varphi$ is a singular functional  there holds
$$
\|x^{\ast}\|=\|v\|_{p'(\cdot)}^{0}+\|\varphi\|,
$$
where $\|v\|_{p'(\cdot)}^{0}$ is the Orlicz norm of $v$ in $L^{p'(\cdot)}([0,1]).$
\end{prop}

\section{Proof of Theorem 1.1}

Let the space $C([0,1])$ is closed subspace in $L^{p(\cdot)}.$ Then there exists a positive constant $c>0$ such that
\begin{equation}
c\leq\|\chi_{(a,b)}\|_{p(\cdot)}\,\,\,\,\mbox{whenever}\,\,\,0\leq a<b\leq1,
\end{equation}
(see \cite{EdGoKo}). It is obvious that for some constant $C>0$
\begin{equation}
\|\chi_{(a,b)}\|_{p(\cdot)}\leq C\,\,\,\,\mbox{whenever}\,\,\,0\leq a<b\leq1.
\end{equation}

From (3.1) and (3.2) we may deduce that for some constants $c_{1},\,c_{2}>0$ and for any $x\in C([0,1])$
\begin{equation}
c_{1}\|x\|_{C}\leq\|x\|_{p(\cdot)}\leq c_{2}\|x\|_{C}
\end{equation}
(see for more details in \cite{EdGoKo}).

Denote $X=C([0,1])$ and $Y=E^{p(\cdot)}([0,1])$ ($X$ is the Banach space with both norms $\|\cdot\|_{C}$ and $\|\cdot\|_{p(\cdot)}$).
 We have  $X\cap Y=\{0\}$ (by (3.1)). Consider the Cartesian product $X\times Y,$ equipped with coordinate-wise vector space operations. For this vector space we have the Banach norm
$$
\|(u,v)\|_{\infty}=\max\{\|u\|_{p(\cdot)},\|v\|_{p(\cdot)}\}.
$$
 Denote the $X\times Y$ vector space equipped with the norm $\|(\cdot,\cdot)\|_{\infty}$ as  $(X\times Y)_{\infty}$. Obviously $(X\times Y)_{\infty}$ is the Banach space. Our main goal is to prove that the mapping  $(X\times Y)_{\infty}\rightarrow X+Y\subset L^{p(\cdot)}([0,1]):$ $(u,v)\rightarrow u+v$ is a (topological) isomorphism. In this case vector space $X+Y$ with the norm $\|\cdot\|_{p(\cdot)}$ is the topological direct sum of the Banach spaces  $X$ and $Y$, and it is written as $X+Y=X\oplus Y.$  From this fact we obtain that the vector space $X+Y=C([0,1])\oplus E^{p(\cdot)}([0,1])$ with the norm $\|\cdot\|_{p(\cdot)}$ is Banach subspace of $L^{p(\cdot)}([0,1])$ and we have
  \begin{equation}
 \|x+y\|_{p(\cdot)}\approx\max\{\|x\|_{p(\cdot)},\|y\|_{p(\cdot)}\},\,\,x\in X,\,y\in Y.
 \end{equation}

 It is obvious that the vector space $X+Y$ with the norm $\|\cdot\|_{p(\cdot)}$ is topological direct sum of Banach spaces $X$ and $Y$ if the  linear projection $P:X+Y\rightarrow X$ which is defined by $P(x+y)=x$ is continuous when $x\in X$ and $y\in Y$. Note that this condition is equivalent to the following: there exists a positive real number $\delta$ such that $\|x-y\|_{p(\cdot)}\geq\delta$ whenever $x\in X, y\in Y$ and $\|x\|_{p(\cdot)}=1.$

Let $x\in X$ and $\|x\|_{p(\cdot)}=1$.  Take  $t_{0}\in[0,1]$  such that
 $$
 |x(t_{0})|=\max_{t\in[0,1]}|x(t)|=\|x\|_{C}.
 $$
  By (3.3) we have
\begin{equation}
1/c_{2}\leq |x(t_{0})|\leq 1/c_{1}.
\end{equation}

We will prove that
$$
d(x)=\inf_{y\in Y}\|x-y\|_{p(\cdot)}\geq \delta>0
$$
for some constant $\delta$ independent of $x$.

By Proposition 2.1 we have
\begin{equation}
d(x)=\lim_{n\rightarrow\infty}\|x-x^{(n)}\|_{p(\cdot)},
\end{equation}
where $x^{(n)}=\chi_{\Omega_{n}}$, $\Omega_{n}=\{t:\,p(t)\leq n\}.$

denote $O_{n}=(t_{0}-\varepsilon_{n},t_{0}+\varepsilon_{n})$, where  numbers $\varepsilon_{n}>0$ we will choose later.

We have
$$
\|x-x^{(n)}\|_{p(\cdot)}=\|x-x^{(n)}\chi_{[0,1]\backslash O_{n}}-x^{(n)}\chi_{[0,1]\cap O_{n}}\|_{p(\cdot)}
$$
\begin{equation}
\geq|\|x-x^{(n)}\chi_{[0,1]\backslash O_{n}}\|_{p(\cdot)}-\|x^{(n)}\chi_{[0,1]\cap O_{n}}\|_{p(\cdot)}|.
\end{equation}

Since for fixed $n$ on the set $\Omega_{n}$  we have $p(t)\leq n$  we may take $O_{n}$ such that $\|x^{(n)}\chi_{[0,1]\cap O_{n}}\|_{p(\cdot)}$ is arbitrary small. Using (3.1) and (3.5) we can choose $O_{n}$ such small that
\begin{equation}
\|x\chi_{O_{n}}\|_{p(\cdot)}\geq \frac{1}{2}|x(t_{0})|\|\chi_{O_{n}}\|_{p(\cdot)}\geq\frac{c}{2c_{2}}.
\end{equation}
From (3.7) and (3.8) we obtain
$$
\|x-x^{(n)}\|_{p(\cdot)}\geq\frac{c}{2c_{2}},
$$
and consequently by (3.6) we have $d(x)\geq\delta=\frac{c}{2c_{2}}.$

 Let $x^{\ast}$ be any  continuous linear  functional from $X^{\ast}.$   It is obvious that $x^{\ast}$ is continuous linear  functional on the space $X$ with the norm $\|\cdot\|_{p(\cdot)}$ (by (3.3)).  Since the space $X\oplus Y$ is a Banach space with the norm $\|\cdot\|_{p(\cdot)}$, the trivial extension (i.e $x^{\ast}(x)=0$, for $x\in Y$) of $x^{\ast}$ is also continuous linear functional on $X\oplus Y$ (see 3.4).  For functional obtained in this way (defined on $X\oplus Y$)
there exist continuous linear  extension (non unique)  on whole space $L^{p(\cdot)}([0,1])$. It is obvious that obtained functional is singular (it is identically $0$ on $E^{p(\cdot)}([0,1]))$ on $L^{p(\cdot)}([0,1])$.
 $\Box$

\textbf{Remark 1}. A closed subspace $Y$ of Banach space $X$ is $M$-ideal in $X$ if $Y^{\bot}$  is the range of bounded projection $P:X^{\ast}\rightarrow X^{\ast}$ which satisfies
$$
\|x^{\ast}\|=\|Px^{\ast}\|+\|x^{\ast}-Px^{\ast}\|\,\,\,\,\mbox{for\,\,all}\,\,\,x^{\ast}\in X^{\ast}.
$$
For more details of general $M$-ideal theory and their applications, we refer to \cite{Ha}. If subspace $Y$ is $M$-ideal in $X$, then $Y$ is proximinal in $X$ (see \cite{Ha}, p.57, Proposition 1.1), that is for any $x\in X$ there exists $y\in Y$ such that
$$
d(x)=\inf_{z\in Y}\|x-z\|=\|x-y\|.
$$

From proposition 2.3 we obtain that the space $E^{p(\cdot)}([0,1])$ is $M$-ideal in $L^{p(\cdot)}([0,1])$ and consequently $E^{p(\cdot)}([0,1])$ is proximinal in $L^{p(\cdot)}([0,1])$; that is
for any $x\in L^{p(\cdot)}([0,1])$ there exists $y\in E^{p(\cdot)}([0,1]) $ such that $d(x)=\|x-y\|_{p(\cdot)}.$

\textbf{Remark 2}.  Let $C([0,1])$ is closed subspace in $L^{p(\cdot)}([0,1])$. Denote  $I=L^{\infty}([0,1])\cap E^{p(\cdot)}([0,1]).$ Note that if $x_{n}\in I, n\in \mathbb{N}$ and $\lim_{n\rightarrow\infty}\|x_{n}-x\|_{\infty}=0,$ then $x\in E^{p(\cdot)}([0,1]).$ It is easy to show that $I$ is an order ideal, which means that it is a closed subspace of $L^{\infty}([0,1])$ with the ideal property.

Note that there exists $\delta>0$ such that  for $x\in C([0,1]), \,\, \|x\|_{C}=1$ and $y\in{I}$ we have $\|x-y\|_{\infty}\geq\delta.$ The last inequality can be proved analogously as it was done in the proof of Theorem 1.1. (It is sufficient to use the inequality $\|x-y\|_\infty \geq\|x-y\|_{p(\cdot)} $ and the fact that $\|x\|_{p(\cdot)}\approx1$). Consequently, vector space $C([0,1])+I$  is the topological direct sum of Banach spaces $C([0,1])$ and $I$ in $L^{\infty}([0,1])$.


\begin{thebibliography}{bif}


\bibitem{AW} Y. A. Abramovich and A. W. Wickstead, Singular extensions and restrictions of order continuous
functionals, Hokkaido Math. J., 21 (1992), 475–482.

\bibitem{AB}  C. D. Aliprantis and O. Burkinshaw, Positive Operators, Pure and Applied  Mathematics, vol. 119, Academic Press, New York, 1985.

\bibitem{ACCUO} A. Amenta, J. Conde-Alonso, D. Cruz-Uribe, and J. Ocsariz, On the dual of variable Lebesgue space with unbounded exponent, arXiv:1909.05987, (September, 2019).

\bibitem{CUF} D.~Cruz-Uribe, A.~Fiorenza,
 Variable Lebesgue Spaces. Foundations and Harmonic Analysis,
Birkh\"auser, Basel (2013).

\bibitem{DHHR} L.~Diening, P.~Harjulehto, P.~H\"ast\"o, M.~R$\stackrel{\circ}{\mbox{u}}$\v zi\v cka,
 Lebesgue and Sobolev spaces with variable exponents, Lecture Notes in Mathematics, 2017. Springer, Heidelberg (2011).

\bibitem{EdGoKo} D. Edmunds, A. Gogatishvili, and T. Kopaliani, Construction of function spaces close to $L^{\infty}$
with associate space close to $L^{1}$, J. Fourier Anal. Appl., 24(6):1539-1553, 2018.

\bibitem{Ha} P. Harmand, D. Werner,and W. Werner, $M$-ideals in Banach Spaces and Banach
Algebras, Lecture Notes in Math. 1547, Springer-Verlag, Berlin-Heidelberg (1993).

\bibitem{He} W. Hensgen, An example concerning the Yosida-Hewitt decomposition of finitely additive measures,
Proc. Amer. Math. Soc. 121 (1994), 641-642.

\bibitem{Hud}  H. Hudzik, Z. Zbjszyniak, Smoothness in Musielak-Orlicz spaces equipped with the Orlicz
norm, Collect. Math., 48 (1997), 543-561.

\bibitem{KoZv} T. Kopaliani and S. Zviadadze, Note on the variable exponent Lebesgue function spaces close
to $L^{\infty}$,  J. Math. Anal. Appl., 474(2):1463-1469, 2019.

\bibitem{Tol} J. F. Toland, Localizing weak convergence in $L_{\infty}$, arXiv:1802.01878 (November 2018).


\bibitem{Val} M. Valadier, Une singuli`ere forme lin\`{e}aire sur $L^{\infty}$. S\'{e}m. D\'{a}nalyse Convexe Montpelier, 4, (1987).

\bibitem{WR} A. J. Wrobel, A sufficient condition for a singular functional on $L^{\infty}[0; 1]$ to be represented on
$C[0; 1]$ by a singular measure, Indag. Math. 29 (2), (2018), 746-751.

\end{thebibliography}
\end{document}